\documentclass[reqno,12pt]{amsart}
\usepackage{amsfonts}
\usepackage{amssymb}
\usepackage{bbm}
\usepackage{times}

\def\B{\mathcal{B}}
\def\D{\mathbb{D}}

\def\C{\mathbb{C}}

\def\B{\mathcal{B}}

\def\msk{\medskip}

\def\bege{\begin{equation}} \def\ende{\end{equation}}
   
\def\b{\beta}   
 
\def\begr{\begin{eqnarray}} \def\endr{\end{eqnarray}}

\def\bege{\begin{equation}} \def\ende{\end{equation}}
\def\begr{\begin{eqnarray}} \def\endr{\end{eqnarray}}
\def\bnum{\begin{enumerate}} \def\enum{\end{enumerate}}

\begin{document}

\title[linear combination of composition operators]{linear combination of composition operators on $H^\infty$ and the Bloch space}
\author{Yecheng Shi and Songxiao Li$^*$ }

\address{Yecheng Shi\\ Faculty of Information Technology, Macau University of Science and Technology, Avenida Wai Long,
Taipa, Macau.
}\email{ 09ycshi@sina.cn}

\address{Songxiao Li\\ Institute of Fundamental and Frontier Sciences, University of Electronic Science and Technology of China,
610054, Chengdu, Sichuan, P. R. China.  \newline
Institute of Systems Engineering, Macau University of Science and Technology, Avenida Wai Long, Taipa, Macau. } \email{jyulsx@163.com}

\begin{abstract} Let $\lambda_i (i=1,...,k)$ be any nonzero complex scalars and $\varphi_i (i=1,..,k)$ be any analytic self-maps of the unit disk $\mathbb{D}$. We show that  the  operator $\sum_{i=1}^k\lambda_iC_{\varphi_i}$ is compact on the Bloch space $\mathcal{B}$ if and only if
 $$\lim_{n\to\infty}\|\lambda_1\varphi_1^n+\lambda_2\varphi_2^n+...+\lambda_k\varphi_k^n\|_{\mathcal{B}}=0.$$
 We also study the linear combination of composition operators on the Banach algebra of bounded analytic functions.

\thanks{*Corresponding author.}
\thanks{This  project was partially supported by NSF of China (No.11471143 and No.11720101003)
   }
\end{abstract}

\maketitle

\section{Introduction}
Let $\D$ denote the open unit disk in the complex plane $\C$. We denote by $H(\D)$ the class of functions analytic in $\D$.
Let $H^\infty = H^\infty(\D)$ be the space of all bounded analytic functions on $\D$.
Then $H^\infty$ is a Banach algebra with the supremum norm $\|f\|_{\infty} = \sup_{z\in\D}  |f(z)| .$

Recall that the Bloch space $\mathcal{B}$ is a space which consists of all $f\in H(\mathbb{D})$ such that
$$
\|f\|_{\b}=\sup_{z \in \mathbb{D}}(1-|z|^2) |f'(z)| < \infty.
$$
It is well known that $\B$ is a Banach space under the norm $ \|f\|_{\B}=|f(0)|+\|f\|_{\b}$.

 For $a\in\D$, let $\sigma_{a}(z):=\frac{a-z}{1-\bar{a}z}$ be the disc automorphism that exchanges 0 for $a$. For $z$, $w\in\D$, the pseudo-hyperbolic distance between $z$ and $w$ is given by
$$\rho(z,w)=|\sigma_{w}(z)|=\bigg|\frac{z-w}{1-\bar{w}z}\bigg|.$$
We write $\rho_{ij}(z)=\rho(\varphi_i(z),\varphi_j(z))$  and denote the hyperbolic derivative $$\varphi^{\#}(z)=\frac{1-|z|^2}{1-|\varphi(z)|^2}\varphi^{\prime}(z).$$
 Let $\bigtriangleup$ denote the collection of all sequences $\{z_n\}$ in
$\D$ converging to some point of $\partial\D$ such that the sequences $\{\varphi_{i}(z_n)\}$, $\{\varphi_i^{\#}(z_n)\}$
and $\{\rho_{ij}(z_n)\}$ also converge for all $i,j=1,...,k$. Given a sequence $\{z_n\}\in \bigtriangleup$ and
an index $j=1,...,k$, define
$$I\{z_n\}=\{i:|\varphi_i(z_n)|\to1\},$$
$$I_j\{z_n\}=\{i:|\rho_{ij}(z_n)|\to0\},$$
$$I^*_j\{z_n\}=I_j\{z_n\}\cap\{i:|\varphi_i^{\#}(z_n)|\nrightarrow 0\},$$
$$I^{\#}_j\{z_n\}=I_j\{z_n\}\cap\{i: \lim\varphi_i^{\#}(z_n)=\lim\varphi_j^{\#}(z_n)\}.$$
Note that every sequence $\{z_n\}$ in $\D$ with $|z_n|\to1$ has subsequences belonging to $\bigtriangleup$.
Furthermore, the sets $I_j\{z_n\}$ induce a natural partition of $I\{z_n\}$.

Let $\varphi$ be an analytic self-map of $\D$. The function $\varphi$ induces a composition operator $C_\varphi:H(\D)\rightarrow H(\D)$   defined by
$C_\varphi f=f\circ\varphi$. An extensive study on the
theory of composition operators has been established during the past four decades. 
 A basic and interesting problem concerning  composition operators is to relate operator theoretic properties to their function theoretic properties of their symbols. We refer the reader to
  \cite{cm} and \cite{zhu}.

 Madigan and Matheson in \cite{mm} proved that $C_\varphi : \B \rightarrow \B$ is compact if and only if $\lim_{|\varphi (z)|\rightarrow1} |\varphi^{\#}(z)| = 0.$ In \cite{WZZ}, Wulan, Zheng and Zhu used the power functions as test functions and obtained a new characterization for the compactness of the   operator $C_\varphi : \B \rightarrow \B$, i.e., they showed that $C_\varphi : \B \rightarrow \B$ is compact if and only if $$\lim_{n\rightarrow\infty}\|\varphi^n\|_{\B}=0.$$
Hosokawa and Ohno \cite{HO} and Nieminen \cite{Ni} have investigated the compact differences of composition operators on the Bloch space. In \cite{SL}, the authors of this paper obtained several estimates for the essential norm of the differences of composition operators on $\B$.  Among others, we showed that
$$\|C_\varphi-C_\psi\|_{\B\rightarrow \B,e}\asymp \lim_{n\to\infty}\|\varphi^n-\psi^n\|_{\B}.$$
For further results on compact differences on various  settings, we refer to \cite{CKM,Mo,Ni,NS,S1,S2, SL} and references therein.




Along the line of study on differences, the study on linear combination of composition operators has been a topic of growing interest.
 Izuchi and Ohno \cite{IO} have given a complete characterization for the compactness of linear combination of composition operators on $H^\infty$.  Hosokawa, Nieminen and Ohno \cite{HNO} characterized the compactness of the linear combination of composition operators on the Bloch space. The main result in \cite{HNO} is stated as follows. \msk

  \noindent{\bf Theorem A \cite{HNO}.}  {\it Suppose $k$ is a positive integer. Let $\lambda_1, \lambda_2, ... , \lambda_k$ be any nonzero complex scalars, and $\varphi_1, \varphi_2,..., \varphi_k$ be any analytic self-maps of $\D$. Then the following statements are equivalent:

  (i) The operator $\sum_{i=1}^k\lambda_iC_{\varphi_i}$ is compact on $\B$.

(ii) $\sum_{i\in I_j\{z_n\}}\lambda_i\varphi_i^{\#}(z_n)\rightarrow 0$ as $n\rightarrow\infty$ for all $\{z_n\}\in\bigtriangleup$, $j\in I\{z_n\}$.

(iii)  $\sum_{i\in I^{*}_j\{z_n\}}\lambda_i= 0$ for all $\{z_n\}\in\bigtriangleup$, $j\in I\{z_n\}$.

(iv) $\sum_{i\in I^{\#}_j\{z_n\}}\lambda_i = 0$ for all $\{z_n\}\in\bigtriangleup$, $j\in I\{z_n\}$ with $\varphi^{\#}_j(z_n)\nrightarrow0$.}

 For further results on the linear combination of composition operators, we refer to \cite{CCK, CIK,CKM,CKWY,GM, HNO, IO,KW2, KM} and references therein.


 The main aim in this work is to provide criteria for the  compactness of the  linear combination of composition operators $\sum_{i=1}^k\lambda_iC_{\varphi_i}$ acting on the Bloch space and $H^\infty$  in terms of the sequence $\|\lambda_1\varphi_1^n+\lambda_2\varphi_2^n+...+\lambda_k\varphi_k^n\|_{\B} $ and $ \|\lambda_1\varphi_1^n+\lambda_2\varphi_2^n+...+\lambda_k\varphi_k^n\|_{\infty}$, respectively. The compactness of the operator $\sum_{i=1}^k\lambda_iC_{\varphi_i}$ on the Bloch space and $H^\infty$ can be easily verified by using the results in this paper.


For two quantities $A$ and $B$, we use the abbreviation
$A\lesssim B$ whenever there is a positive constant $c$ (independent of the associated variables) such that $A\leq cB$.
We write $A\asymp B$, if $A\lesssim B\lesssim A$.


\section{ Main results and proofs}
In this section, we state and prove some preliminary results needed for the rest of this paper. The following  lemma below is an adaptation of in \cite{WZZ} to a more general and abstract setting.\msk

  \noindent{\bf  Lemma 1.}  {\it Let $(X,\|\cdot\|_{X})$ and $(Y,\|\cdot\|_{Y})$ be two Banach spaces of analytic functions in $\D$ and
$\|\cdot\|_{X\rightarrow Y}$ is the operator norm. $T$ is a bounded linear operator from $X$ into $Y$. Suppose $\{f_n\}$ is a function sequence of $X$,  $f_n(z)=\sum_{j=0}^\infty a_{n,j}z^j$, such that

  (i) there exists a positive constant $C$ independent of $n$ such that
  $$\sum_{j=0}^\infty |a_{n,j}|\|z^j\|_{X}\leq C,$$

  (ii) for any fix positive integer $k$,
  $$\lim_{n\to\infty}\sum_{j=0}^k |a_{n,j}|\|z^j\|_{X}=0.$$
Then
$$\lim_{n \to \infty} \|Tf_n\|_{Y} \lesssim \limsup_{n\to\infty}\frac{\|Tz^n\|_{Y}}{\|z^n\|_{X}}.$$}

\noindent{\it Proof.} Since $f_n(z)=\sum_{j=0}^\infty a_{n,j}z^j$, we have
\begr
\|Tf_n\|_Y &\leq& \sum_{j=0}^\infty |a_{n,j}|\|Tz^j\|_Y\nonumber\\
&=& \sum_{j=0}^\infty |a_{n,j}|\|z^j\|_X\frac{\|Tz^j\|_Y}{\|z^j\|_X}\nonumber\\
&\leq& \sum_{j=0}^k |a_{n,j}|\|z^j\|_X\frac{\|Tz^j\|_Y}{\|z^j\|_X}+\big(\sum_{j=k+1}^\infty |a_{n,j}|\|z^j\|_X\big)\sup_{i\geq k+1}\frac{\|Tz^i\|_Y}{\|z^i\|_X}\nonumber\\
&\leq& \big(\sum_{j=0}^k |a_{n,j}|\|z^j\|_X\big)\|T\|_{X\to Y}+\big(\sum_{j=k+1}^\infty |a_{n,j}|\|z^j\|_X\big)\sup_{i\geq k+1}\frac{\|Tz^i\|_Y}{\|z^i\|_X}.\nonumber
\endr
By  $\it (i)$ and $\it (ii)$, letting $n\to\infty$, we have
 $$\limsup_{n\to\infty}\|Tf_n\|_Y \lesssim \sup_{i\geq k+1}\frac{\|Tz^i\|_Y}{\|z^i\|_X}.$$
Letting $k\to\infty$, we get the desired result. The proof is complete.\msk


\noindent{\bf Lemma 2. }{\it  Let $f, g, h \in H(\D)$. Suppose $f(z)=\sum_{i=0}^\infty a_iz^i$,
 $g(z)=\sum_{j=0}^\infty b_jz^j$, and  $h(z)=\sum_{k=0}^\infty c_kz^k$, where $a_i, b_j, c_k \in \C$. If $h(z)=f(z)g(z)$, then
 $$\sum_{k=0}^\infty|c_k|\leq \Big(\sum_{i=0}^\infty |a_i|\Big)\Big(\sum_{i=0}^\infty|b_i|\Big).$$
 }
\noindent{\it Proof.} Since
\begr
h(z)=\Big(\sum_{i=0}^\infty a_iz^i\Big)\Big(\sum_{j=0}^\infty b_jz^j\Big)
=\sum_{k=0}^\infty\Big(\sum_{i=0}^k a_ib_{k-i}\Big)z^k,\nonumber
\endr
we obtain
$c_k=\sum_{i=0}^ka_ib_{k-i}.$ Thus,
\begr
\sum_{k=0}^\infty|c_k|&=&\sum_{k=0}^\infty\Big|\sum_{i=0}^ka_ib_{k-i}\Big|
\leq \sum_{k=0}^\infty\sum_{i=0}^k|a_i||b_{k-i}|
= \Big(\sum_{i=0}^\infty |a_i|\Big)\Big(\sum_{i=0}^\infty|b_i|\Big).\nonumber
\endr\msk


 \noindent{\bf Theorem 1.}  {\it Suppose $k$ is a positive integer. Let $\lambda_1, \lambda_2, ... , \lambda_k$ be any nonzero complex scalars, and $\varphi_1, \varphi_2,..., \varphi_k$ be any analytic self-maps of $\D$. Then the  operator
$\sum_{i=1}^k\lambda_iC_{\varphi_i}$ is compact on $\B$ if and only if
\begr \lim_{n\to\infty}\|\lambda_1\varphi_1^n+\lambda_2\varphi_2^n+...+\lambda_k\varphi_k^n\|_{\B}=0.\endr}

\noindent{\it Proof.}  Suppose $\sum_{i=1}^k\lambda_iC_{\varphi_i}$ is compact on $\B$. Consider the test functions $p_n(z)=z^n$. We have
$\|p_{n}\|_\B\leq\frac{2}{e}$ and $f_n\rightarrow0$ weakly in $\B$ as $n\rightarrow\infty$.  Then by the compactness of $\sum_{i=1}^k\lambda_iC_{\varphi_i}$, we get (1), as desired.

Conversely, suppose that (1) holds.  By Theorem A, $\sum_{i=1}^k\lambda_iC_{\varphi_i}$ is compact on $\B$ if and only if
$$\sum_{i\in I_j\{z_n\}}\lambda_i\varphi_i^{\#}(z_n)\to0 \mbox{~as~}n\to\infty \mbox{~for all~}\{z_n\}\in\bigtriangleup, j\in I\{z_n\}.$$
Suppose that $\{z_n\}\in \bigtriangleup$ such that $|\varphi_j(z_n)|\to1$, we write $I=I\{z_n\}$ and $J=I_{j}\{z_n\}$.
Let $\{f_n\}$ be the sequence of analytic functions defined by
$$f_n(z)=\sigma_{\varphi_{j}(z_n)}(z)\prod_{i\in I\backslash J}\sigma_{\varphi_i(z_n)}(z)^2-\gamma_n,$$
where $\gamma_n=\varphi_{j}(z_n)\prod_{i\in I\backslash J}\varphi_i(z_n)^2$.
By the proof of Theorem A (see \cite{HNO}), we see that
if $\lim_{n\to\infty}\|\sum_{i=1}^k\lambda_iC_{\varphi_i} f_n\|_{\B}=0$, then
$$\lim_{n\to\infty}\sum_{i\in J}\lambda_i\varphi_i^{\#}(z_n)=0.$$
We denote $a_{n,i} = \varphi_{i}(z_n)$, for $i\in I$. Then
$a_{n,i}\to1$ as $n\to\infty$ and
$$f_n(z):=\sigma_{a_{n,j}}(z)\prod_{i\in I\backslash J}\sigma_{a_{n,i}}(z)^2-\gamma_n, \mbox{~~and~~} \gamma_n = a_{n,j}\prod_{i\in I\backslash J}a^2_{n,i}.$$
The sufficiency of the proof will be given after we prove that
$$\lim_{n\to\infty}\|\sum_{i=1}^k\lambda_iC_{\varphi_i} f_n\|_{\B} \lesssim \limsup_{n\to\infty}\|\lambda_1\varphi_1^n+...\lambda_k\varphi_k^n\|_{\B}.$$
Suppose 
  $f_n(z)=\sum_{l=1}^\infty b_{n,l}z^l.$ Since $\sum_{i=1}^k\lambda_iC_{\varphi_i} : \B \to \B$ is bounded and $\|z^j\|_{\B}\leq \frac{2}{e}$, by Lemma 1, it is enough to prove that for any fix positive integral $N$,
$\lim_{n\to\infty}\sum_{l=1}^N |b_{n,l}| = 0$
and there exists a positive constant $C$ independent of $n$ such that
$\sum_{l=1}^\infty |b_{n,l}|<C.$

Since $\sigma_{a}(z)=a-(1-|a|^2)\sum_{l=0}^{\infty}\bar{a}^lz^{l+1},$
we get that
\begr
&&\big(a_{n,j}-(1-|a_{n,j}|^2)\sum_{l=0}^{N-1}{\overline{a_{n,j}}}^lz^{l+1}\big) \prod_{i\in I\backslash J}\Big(a_{n,i}-(1-|a_{n,i}|^2)\sum_{l=0}^{N-1}{\overline{a_{n,i}}}^lz^{l+1}\Big)^2\nonumber\\
&=&\gamma_n+b_{n,1}z^1+b_{n,2}z^2+...+b_{n,N}z^{N}+z^{N+1}h(z),\nonumber
\endr
where $h(z)$ is a polynomial.

Using Lemma 2, we obtain
\begr
&&\sum_{l=1}^{N}|b_{n,l}|\nonumber\\
&\leq& \big(|a_{n,j}|+(1-|a_{n,j}|^2)\sum_{l=0}^{N-1}{|a_{n,j}|}^l\big) \prod_{i\in I\backslash J}\Big(|a_{n,i}|+(1-|a_{n,i}|^2)\sum_{l=0}^{N-1}{|a_{n,i}|}^l\Big)^2\nonumber\\
&~~~~~~~~~~~~~~&-|\gamma_n|\to0 \mbox{~~as~~} n\to \infty \nonumber
\endr
and
\begr
&&\sum_{l=1}^{\infty}|b_{n,l}|\leq\sum_{l=1}^{\infty}|b_{n,l}|+|\gamma_n|\nonumber\\
&\leq& \big(|a_{n,j}|+(1-|a_{n,j}|^2)\sum_{l=0}^{\infty}{|a_{n,j}|}^l\big) \prod_{i\in I\backslash J}\Big(|a_{n,i}|+(1-|a_{n,i}|^2)\sum_{l=0}^{\infty}{|a_{n,i}|}^l\Big)^2\nonumber\\
&\leq& 3^{2|I\backslash J|+1}.
 \nonumber
\endr
Here we used the fact that for all $a\in\D$, $|a|+(1-|a|^2)\sum_{l=0}^{\infty}{|a|}^l\leq3$ and  $|X|$ denote the cardinal number of $X$.
The proof is complete.\msk

Combining with Corollary 3.3 of \cite{HNO}, we have the following results.\msk

\noindent{\bf Corollary 1.}  {\it Let $\lambda_1, \lambda_2$ be nonzero complex scalars, and $\varphi, \psi$ be any two analytic self-maps of $\D$. Suppose that none of $C_\varphi$ and $C_\psi$ is compact on $\B$. Then the  operator $\lambda_1C_{\varphi}+\lambda_2C_{\psi}$ is compact on $\B$ if and only if $\lambda_1+\lambda_2=0$
and $\lim_{n\to\infty}\|\varphi^n-\psi^n\|_{\B}=0.$}\msk

\noindent{\bf Corollary 2.}  {\it Let $\lambda_1, \lambda_2$ be nonzero complex scalars such that $\lambda_1+\lambda_2\neq 0$. Let $\varphi, \psi$ be any two analytic self-maps of $\D$. Then
$\lim_{n\to\infty}\|\lambda_1\varphi^n+\lambda_2\psi^n\|_{\B}=0$ if and only if    $ \lim_{n\to\infty}\|\varphi^n\|_{\B} =0$ and $ \lim_{n\to\infty}\|\psi^n\|_{\B} =0.$     }\msk

 \noindent{\bf Theorem 2.}  {\it Suppose $k$ is a positive integer. Let $\lambda_1, \lambda_2, ... , \lambda_k$ be any nonzero complex scalars, and $\varphi_1, \varphi_2,..., \varphi_k$ be any analytic self-maps of $\D$. Then the  operator $\sum_{i=1}^k\lambda_iC_{\varphi_i}$ is compact on $H^\infty$ if and only if
\begr \lim_{n\to\infty}\|\lambda_1\varphi_1^n+\lambda_2\varphi_2^n+...+\lambda_k\varphi_k^n\|_{\infty}=0.\endr}

\noindent{\it Proof.}  It is easy to see that
$$\|\sum_{i=1}^k\lambda_iC_{\varphi_i}\|_{H^\infty\to H^\infty}\leq\sum_{i=1}^k|\lambda_i|\|C_{\varphi_i}\|_{H^\infty\to H^\infty}\leq \sum_{i=1}^k|\lambda_i|.$$
Suppose $\sum_{i=1}^k\lambda_iC_{\varphi_i}$ is compact on $H^\infty$.
Consider the test functions $p_n(z)=z^n$. We have
$\|p_{n}\|_{\infty} = 1$ and $p_n\rightarrow 0$ uniformly on any compact subsets of $\D$ as $n\to\infty$. Then (see Proposition 2.1 of \cite{IO}),
 (2) holds.

Conversely, assume that (2) holds.  By the proof of Theorem 2.2 of \cite{IO}, $\sum_{i=1}^k\lambda_iC_{\varphi_i}$ is compact on $H^\infty$ if and only if
$$\lim_{n\to\infty}\|\sum_{i=1}^k\lambda_iC_{\varphi_i} g_n\|_{\infty}=0 \mbox{~as~}n\to\infty \mbox{~for all~}\{z_n\}\in\bigtriangleup, j\in I\{z_n\},$$
where
$$g_n(z)=\frac{1-|\varphi_{j}(z_n)|^2}{1-\overline{\varphi_j(z_n)}z}\prod_{i\in I\backslash J}\sigma_{\varphi_i(z_n)}(z).$$
Here  $I=I\{z_n\}$ and $J=I_{j}\{z_n\}$. We denote $a_{n,i} = \varphi_{i}(z_n)$, for all $i\in I$. Then
$$g_n(z)=\frac{1-|a_{n,j}|^2}{1-\overline{a_{n,j}}z}\prod_{i\in I\backslash J}\sigma_{a_{n,i}}(z).$$
Suppose $g_n(z)=\sum_{l=0}^\infty c_{n,l}z^l.$  Then
\begr
&&\big((1-|a_{n,j}|^2)\sum_{l=0}^{N}{\overline{a_{n,j}}}^lz^{l}\big) \prod_{i\in I\backslash J}\Big(a_{n,i}-(1-|a_{n,i}|^2)\sum_{l=0}^{N-1}{\overline{a_{n,i}}}^lz^{l+1}\Big)\nonumber\\
&=&c_{n,0}+c_{n,1}z^1+c_{n,2}z^2+...+c_{n,N}z^{N}+z^{N+1}p(z),\nonumber
\endr
where $p(z)$ is a polynomial.

Using Lemma 2 again, we obtain
\begr
&&\sum_{l=0}^{N}|c_{n,l}|\nonumber\\
&\leq& \big((1-|a_{n,j}|^2)\sum_{l=0}^{N}{|a_{n,j}|}^l\big) \prod_{i\in I\backslash J}\Big(|a_{n,i}|+(1-|a_{n,i}|^2)\sum_{l=0}^{N-1}{|a_{n,i}|}^l\Big)\to0\nonumber\\
&~~~~~~~~~~~~~~& \mbox{~~as~~} n\to \infty \nonumber
\endr
and
\begr
\sum_{l=0}^{\infty}|c_{n,l}|
&\leq& \big((1-|a_{n,j}|^2)\sum_{l=0}^{\infty}{|a_{n,j}|}^l\big) \prod_{i\in I\backslash J}\Big(|a_{n,i}|+(1-|a_{n,i}|^2)\sum_{l=0}^{\infty}{|a_{n,i}|}^l\Big)\nonumber\\
&=& (1+|a_{n,j}|) \prod_{i\in I\backslash J}\big(|a_{n,i}|+(1+|a_{n,i}|)\big)\nonumber\\
&\leq&2\cdot3^{|I\backslash J|}. \nonumber
\endr
Using Lemma 1, we get
$$\lim_{n\to\infty} \|\sum_{i=1}^k\lambda_iC_{\varphi_i} g_n\|_{\infty}\lesssim \limsup_{n\to\infty}\|\lambda_1\varphi_1^n+...\lambda_k\varphi_k^n\|_{\infty}=0,$$
as desired. The proof is complete.\\





Next, we give some further remarks for Theorem A and Theorem 2.2 of \cite{IO}.
Set $\Lambda_k:=\{1,2,..,k\}$.\msk

\noindent{\bf Theorem 3.}  {\it Suppose $k$ is a positive integer. Let $\varphi_1, \varphi_2,..., \varphi_k$ be analytic self-maps of $\D$ and $\lambda_1, \lambda_2, ... , \lambda_k$ be nonzero complex scalars. Suppose  $\sum_{i\in J}\lambda_i\neq 0$ for every non-empty proper subset $J$ of $\Lambda_k$  and at least one of the operator $C_{\varphi_1},...,C_{\varphi_k}$ is not compact on $\B$. Then the operator $\sum_{i=1}^k\lambda_iC_{\varphi_i}$ is compact on $\B$ if and only if the following two conditions are satisfied:

(i) $\sum_{j=1}^k \lambda_j=0$;

(ii) $C_{\varphi_i}-C_{\varphi_j}$ is compact for each $i,j\in\Lambda_k$. }\msk

\noindent{\it Proof.} Suppose that $\sum_{i=1}^k\lambda_iC_{\varphi_i}$ is compact on $\B$. We first prove that for each $\{z_n\}\in\bigtriangleup$
$$I_1^{*}\{z_n\}=...=I_k^{*}\{z_n\}=\emptyset \mbox{~~or~~} I_1^{*}\{z_n\}=...=I_k^{*}\{z_n\}=\Lambda_k.$$
In fact, for each $\{z_n\}\in \bigtriangleup$, there exists $i\in\Lambda_k$ such that $I_i^{*}\{z_n\}$ is nonempty.
Then by Theorem A $\it(iii)$, we get $\sum_{i\in I^{*}_i\{z_n\}}\lambda_i=0$. Therefore,
  $$I\{z_n\}=I_i\{z_n\}= I_i^{*}\{z_n\}=\Lambda_k.$$
Then $I_1^{*}\{z_n\}=...=I_k^{*}\{z_n\}=\Lambda_k$. Hence,
$$I_1^{*}\{z_n\}=...=I_k^{*}\{z_n\}=\emptyset \mbox{~~or~~} I_1^{*}\{z_n\}=...=I_k^{*}\{z_n\}=\Lambda_k.$$
Therefore,
$$\lim_{|\varphi_i(z)|\to1}\varphi_i^{\#}(z)\rho(
\varphi_i(z),\varphi_j(z))=0$$
and
$$\lim_{|\varphi_j(z)|\to1}\varphi_i^{\#}(z)\rho(
\varphi_i(z),\varphi_j(z))=0$$
for each $i,j\in\Lambda_k$.
By the result in \cite{Ni},  $C_{\varphi_i}-C_{\varphi_j}$ is compact for each $i,j\in\Lambda_k$.

 By the assumption, we assume that $C_{\varphi_j}$ is not compact on $\B$. Then there exists $\{z_n\}\in \bigtriangleup$ such that $I_j^{*}\{z_n\}$ is nonempty.
Then $I^{*}_j\{z_n\}=\Lambda_k$ and
$\sum_{i\in I^{*}_j\{z_n\}}\lambda_i=0$. Therefore,
$C_{\varphi_1},...,C_{\varphi_k}$ are all not compact on $\B$ and
$\sum_{j=1}^k \lambda_j=0.$

Conversely, assume that $\it (i)$ and $\it (ii)$ hold.  Since
$$\sum_{i=1}^k\lambda_iC_{\varphi_i}=(\sum_{i=1}^k\lambda_i)C_{\varphi_1}+\sum_{i=2}^k\lambda_i(C_{\varphi_i}-C_{\varphi_1})=\sum_{i=2}^k\lambda_i(C_{\varphi_i}-C_{\varphi_1}),$$
we have that
$\sum_{i=1}^k\lambda_iC_{\varphi_i}$ is compact on $\B$.\msk

\noindent{\bf Remark 1.}   Suppose $\sum_{i\in J}\lambda_i\neq 0$ for every non-empty proper subset $J$ of $\Lambda_k$. Then $\sum_{i=1}^k\lambda_iC_{\varphi_i}$  is compact on $\B$ if and  only if all of $C_{\varphi_1},...,C_{\varphi_k}$ are compact on $\B$ or none of them is compact and conditions $\it (i)$ and $\it (ii)$  in Theorem 3 hold. \msk

\noindent{\bf Example 1.}   Applying Theorem 3, we can easily check whether a linear combination of composition operators is compact on $\B$ under some mild conditions. For example:

$\bullet~$ $6C_{\varphi_1}-C_{\varphi_2}-2C_{\varphi_3}-3C_{\varphi_4}$ is compact on $\B$ if and only if $C_{\varphi_1}-C_{\varphi_i}$ is compact on $\B$ for each $i=2,3,4.$

$\bullet~$ $4C_{\varphi_1}-C_{\varphi_2}-2C_{\varphi_3}$ is compact on $\B$ only when $C_{\varphi_1}, C_{\varphi_2}, C_{\varphi_3}$ are all compact on $\B$.

$\bullet~$ $3C_{\varphi_1}-iC_{\varphi_2}-2C_{\varphi_3}$ is compact on $\B$ only when $C_{\varphi_1}, C_{\varphi_2}, C_{\varphi_3}$ are all compact on $\B$.\msk

\end{document}